\documentclass[final,1p,times]{elsarticle}

\usepackage{amssymb}
\usepackage{amsmath}
\usepackage{graphicx}
\usepackage{xspace}
\usepackage{soul}  % ten pakiet pozwala uzywac komendy \st{} do przekreslania tekstu
\usepackage[ruled]{algorithm2e}

\journal{Communications in Nonlinear Science and Numerical Simulation}

\newcommand{\capd}{CAPD::DynSys\xspace}
\newcommand{\PM}{\mathcal{P}}

\usepackage{xcolor,listings}
\newtheorem{theorem}{Theorem}

\newenvironment{proof}{\textbf{Proof:}}{\qed}

\begin{document}
\lstset{language=C++,
    keywordstyle=\color{blue}\bfseries,
    commentstyle=\color{gray},
    stringstyle=\ttfamily\color{red!50!brown},
    showstringspaces=false}

\begin{frontmatter}
\title{\capd: a flexible C++ toolbox for rigorous numerical analysis of dynamical systems}

\author{Tomasz Kapela\fnref{MMMaestro}}
\ead{Tomasz.Kapela@uj.edu.pl}

\author{Marian Mrozek\fnref{MMMaestro}}
\ead{Marian.Mrozek@ii.uj.edu.pl}

\author{Daniel Wilczak\fnref{MMMaestro}\corref{ca}}
\ead{Daniel.Wilczak@ii.uj.edu.pl}

\author{Piotr Zgliczy\'nski\fnref{PZ}\corref{}}
\ead{Piotr.Zgliczynski@ii.uj.edu.pl}

\address{Faculty of Mathematics and Computer Science,
Jagiellonian University, \L ojasiewicza 6, 30-348 Krak\'ow, Poland.}
\fntext[MMMaestro] {This research is partially supported by
    the Polish National Science Center under Maestro Grant No. 2014/14/A/ST1/00453 and under Grant No. 2015/19/B/ST1/01454.}
\fntext[PZ] {This research is partially supported by
    the Polish National Science Center under Grant No. 2019/35/B/ST1/00655.}

\cortext[ca] {Corresponding author.}
\date{\today}
\begin{abstract}
We present the \capd library for rigorous numerical analysis of dynamical systems. The basic interface is described together with several interesting case studies illustrating how it can be used for computer-assisted proofs in dynamics of ODEs.
\end{abstract}

\begin{keyword}
rigorous numerical analysis\sep
C++ library\sep
computer-assisted proof


\sep % NA.Algorithms with automatic result verification
 % Periodic orbits of vector fields and flows
%code \sep code (2000 is the default)
\end{keyword}
\end{frontmatter}

\section{Introduction}

In the study of nonlinear ODEs, there is a huge gap between  what we can observe in  numerical simulations and what we can prove rigorously. It is possible to overcome this problem by means of computer-assisted poofs. For its realization it is desirable to have a library for rigorous integration of ODEs and computation of
Poincar\'e maps derived from ODEs. There are several libraries designed for rigorous integration of ODEs. Some of them are open source, just to mention \cite{Nedialkov2006,KV,RauhBrillGunter2009,bresolin-hscc2020}, and some are not \cite{BerzMakino1999}. To the best of our knowledge none of them directly supports computation of Poincar\'e maps, which is a powerful tool for studying dynamics of ODEs.

In the present paper we describe the \capd library \cite{CAPD} which is well suited for this task. What this library offers may be described as follows.

Consider an initial value problem for an  ODE
\begin{eqnarray}
  x'&=&f(\lambda,t,x),  \label{eq:ODE} \\
  x(t_0)&=&x_0, \label{eq:IC-ODE}
\end{eqnarray}
where $x \in \mathbb{R}^n$, $t$ is a time variable, $\lambda \in \Lambda \subset \mathbb{R}^k$ is a fixed parameter and $f:\Lambda \times \mathbb{R} \times \mathbb{R}^n \to \mathbb{R}^n $ is a smooth, 'programmable' function. Let $\varphi(t,t_0,\lambda,x_0)$ be a solution of (\ref{eq:ODE})--(\ref{eq:IC-ODE}) for some fixed $\lambda$. Given a set of initial conditions $Z \subset \mathbb{R}^n$, parameters $\Delta \subset \Lambda$ and a $t > t_0$ we want to:
\begin{itemize}
\item establish that for all $x_0 \in Z$ and $\lambda \in \Delta$ the solution $\varphi(t,t_0,\lambda,x_0)$ is defined,
\item  give a rigorous bound for $D_a \varphi(t,t_0,\lambda,x_0)$ valid for  all $x_0 \in Z$ and $\lambda \in \Delta$,
  where $D_a$ is   the partial derivative operator with respect $x_0$ and/or $\lambda$ of order $r$. The case $r=0$ means that we compute rigorous bounds for $\varphi(t,t_0,\Delta,Z)$.
\end{itemize}
Analogous questions can be asked for Poincar\'e maps for ODEs.

 The \capd library can accomplish these tasks for many interesting systems for $n$ not too large (an example of a `large' value of $n$ up to which the library performed well is $n=80$, which was used for choreographies in $N$-body problem \cite{KapelaZgliczynski2003,KapelaSimo2007}) and the order of the partial derivatives is also not too large, say $r=2,3$ (the library handled $r=5$ in validation of KAM tori \cite{WilczakBarrio2017}) if the sizes of initial conditions are not too big and  the integration time $t-t_0$ is not very large.

The fact that we can compute (enclose) $\varphi(t,t_0,\Delta,Z)$ in single computation is essential for  the computer
assisted proofs, as very often   abstract theorems in dynamics involve  assumptions on the behaviour of solutions on sets, while the single trajectory computations usually do not lead to interesting rigorous statements about the dynamics of the underlining ODE.

 In the following  discussion we will say that we performed $\mathcal{C}^r$ computations if the partial derivatives  up to order $r$ have been computed.

%Let us briefly outline the history of the CAPD.
\subsection{Short history}

CAPD is an acronym for ``Computer Assisted Proofs in Dynamics''. The library was initiated in early 1990's by Marian Mrozek as the tool for the computer assisted proof of chaotic dynamics in the Lorenz system \cite{MischaikowMrozek1995,MischaikowMrozek1998,MischaikowMrozekSzymczak2001}. The present version of the library is split into \capd and CAPD::RedHom \cite{JudaMrozek2014} parts, devoted to the dynamical systems and topology, respectively. The present article focuses on tools from \capd. The DynSys part was developed by Mrozek's Ph. D. students and their descendants at the Jagiellonian University in Krakow, Poland. The most important contributors are (listed more or less chronologically): P. Zgliczy{\'n}ski, P. Pilarczyk, D. Wilczak, T. Kapela.

Papers of Mischaikow, Mrozek and Szymczak on the Lorenz attractor \cite{MischaikowMrozek1995,MischaikowMrozek1998,MischaikowMrozekSzymczak2001} used $\mathcal{C}^0$ computations, only. Other results using $\mathcal C^0$ computations in CAPD library from these early stages of development are
\begin{itemize}
\item symbolic dynamics for the R\"ossler system \cite{Zgliczynski1997},
\item connecting orbits in the Michelson system \cite{Zelawski1999},
\item periodic orbits in the R\"ossler system \cite{Pilarczyk1999}.
\end{itemize}
The early version of $\mathcal C^0$-integrator was based on the logarithmic norms and it  was  slow and inefficient, when compared to the algorithms currently used by the library.

Around 2000 the Lohner algorithm \cite{Lohner1992} (for $\mathcal{C}^0$-computations) and $\mathcal{C}^1$-Lohner type algorithm \cite{Zgliczynski2002}  for efficient $\mathcal C^1$-computations were implemented in the CAPD library. The implementation was designed for very limited types of vector fields, namely degree two polynomials, which include well known R\"ossler \cite{Rossler1976,Rossler79}, Lorenz \cite{Lorenz1963} and Michelson \cite{Mi} systems. Around that time Daniel Wilczak joined the project and implemented these algorithms for general 'programmable' vector fields using the automatic differentiation \cite{RallCorliss1996}.  The  CAPD library was published online for the first time in 2004. Around the year 2008 the $\mathcal C^r$-Lohner algorithm \cite{WilczakZgliczynski2011} was added to the library, and soon after this a  rigorous solver of differential inclusions \cite{KapelaZgliczynski2009} and support for  computation in high precision were written by Tomasz Kapela.

\subsection{Some computer assisted proofs using CAPD library}

The quality of the bounds provided by the \capd library can be judged by looking at the list of computer-assisted proofs in dynamics of ODEs in which it was used.
The list below is incomplete, we focus only on a number of selected applications. The results using $\mathcal{C}^0$- and $\mathcal{C}^1$-computations  include the questions of the existence of periodic orbits and their local uniqueness, the existence of symbolic dynamics, the existence of hyperbolic
invariants sets, the existence of homo- and heteroclinic orbits. Here are some examples
\begin{itemize}
\item symbolic dynamics in the H\'enon-Heiles Hamiltonian \cite{ArioliZgliczynski2001},
\item symbolic dynamics and symmetric periodic orbits in Michelson system \cite{Wi1},
\item homoclinic and heteroclinic connections between Lyapunov orbits and symbolic dynamics in the planar circular restricted three body problem \cite{WilczakZgliczynski2003,WilczakZgliczynski2005},
\item Shilnikov orbits and Bykov cycles in the Michelson system \cite{Wilczak2006},
\item existence of choreographies in Newtonian $N$-body problem \cite{KapelaZgliczynski2003,KapelaSimo2007},
\item hyperbolic Smale-Williams attractor for Kuznetsov System \cite{Wilczak2010},
\item invariant manifolds in the restricted three body problem by Capi\'nski and his coworkers \cite{Capinski2012,CapinskiWasieczko2015},
\item Birkhoff regions of instability in the three body problem, using Aubry-Mather theory, by Galante and Kaloshin \cite{GalanteKaloshin},
\item existence of double spiral attractor in the Chua's circuits by Galias and Tucker \cite{GaliasTucker2019},
\item counting of periodic orbits of flows by Galias and Tucker \cite{GaliasTucker2008,GaliasTucker2011},
\item stability of $N$-body motions forming platonic polyhedra by Fenucci and Gronchi \cite{Fenucci_2018},
%\item Bartha and Tucker work on the fixed points in the Kuramoto-Sivashinsky equation \cite{BARTHA2015339},  \textbf{PZ: nie wiem kto to dodal, ale to nie jest przyklad zastosowania CAPD - to jest self-consistent bounds - trzeba usnac}
\item study of periodic by orbits by Miyaji and Okamoto \cite{MiyajiOkamoto2014},
\item existence of unimodal solutions in the Proudman--Johnson equation by Miyaji and Okamoto \cite{MiyajiOkamoto2019},
\item study of singularities in dynamical systems by Matsue \cite{Matsue2019},
\item applications to rigorous estimates of reachable sets in the context of control theory and robotic by Jaulin and his coauthors \cite{ROHOU2018379,ROHOU201776} and Cyranka et al. \cite{cyranka2017lagrangian},
\item heteroclinic connections in Ohta--Kawasaki Model by Cyranka and Wanner \cite{CyrankaWanner},
\item attracting invariant tori by Capi\'nski, Fleurantin and James \cite{CapinskiFleurantinJames2020}.
\end{itemize}

To address other phenomena, such as bifurcations of periodic
orbits, invariant tori through the KAM theory, nonlinear stability of elliptic periodic orbits, KAM stability etc.
one needs the knowledge of partial derivatives with respect to the initial conditions of the higher order. Using  algorithms from \capd library for $\mathcal{C}^r$-computations the following results
has been obtained
\begin{itemize}
\item global and local bifurcations of periodic orbits and invariant manifolds \cite{KokubuWilczakZgliczynski2007,WilczakZgliczynski2009focm,WalawskaWilczak2019,WilczakZgliczynski2009siads},
\item Arnold diffusion in the restricted three body problem by Capi\'nski and Gidea \cite{CapinskiGidea2020},
\item non-linear stability of elliptic periodic solutions \cite{KapelaSimo2017,WilczakBarrio2017,BarrioWilczak2020},
\item normally hyperbolic invariant manifolds and computer assisted Melnikov method  \cite{CapinskiRoldan2012,CapinskiZgliczynski2017,CAPINSKI20183988}.
\end{itemize}

\subsection{Outline of the paper}

In Section~\ref{sec:interface} we describe the basic interface to \capd library. In Sections~\ref{sec:c0}, \ref{sec:c1} and \ref{sec:cr} we present a list of case studies on how the \capd library can be used in various contexts. Examples are grouped by the maximal order of space derivatives involved -- we refer to them as  $\mathcal C^0$, $\mathcal C^1$ and $\mathcal C^r$ computations. The examples selected here are on the one hand very short (so that it is possible to write out the full C++ code) but on the other side are non-trivial and present diverse spectrum of mathematical problems, where the \capd library may be helpful.

\section{The CAPD library: interface and basic usage}
\label{sec:interface}

The \capd library provides data structures and algorithms designed for analysis of discrete and continuous dynamical systems in finite and infinite dimension. They are written in the spirit of generic programming with high level of abstraction allowing the user to tune or adapt some subroutines for specific problems.

The \capd library provides algorithms for both non-rigorous and rigorous computation. The non-rigorous ones are mainly used for simulation, prototyping or finding approximations of objects we are interested in. They are based on double precision floating point numbers supported by hardware so they are fast but prone to errors coming from rounding, significant bits cancellation, inaccuracy of numerical method etc. On the other hand, rigorous methods are group of algorithms, which compute an outer bounds of the objects we are interested in (like values and derivatives of maps, solutions to IVPs).

Most often used interface of the \capd library is available via the following two header files
\begin{lstlisting}
#include "capd/capdlib.h"   // CAPD library header
#include "capd/mpcapdlib.h" // Multi-precision CAPD header
using namespace capd;
\end{lstlisting}
All types and algorithms are defined in the main namespace \texttt{capd}.

\subsection{Basic arithmetic types and naming convention.}

 In the \capd library the special type \texttt{interval} provides interval arithmetic (see \cite{Moore1966,Tucker11}) and is a base for all data types used in rigorous computations. The precision provided by built-in floating point types is sometimes not sufficient and causes huge overestimation in rigorous computations. The \capd library defines MpFloat and MpInterval types that provide floating point numbers and intervals, respectively, of arbitrary precision. The implementation is based on the MPFR library \cite{MPFR}. The following example shows the basic usage of the above four arithmetic types.

\begin{lstlisting}
#include <iostream>         // C++ standard output libary
#include "capd/capdlib.h"   // CAPD library header 
#include "capd/mpcapdlib.h" // Multi-precision CAPD header
using namespace std;
using namespace capd;       

template <typename T>
T f(T x){
  return sqr(sin(x))*exp(x) + 2*x*(cos(x));
}
int main(){
  cout.precision(17);
    // Non-rigorous computations
  double x = 0.75; // representable number
  cout << "f=" << f(x) << endl;
    // Rigorous computations using
    // interval arithmetics with hardware support (53 mantissa bits) 
  interval iy = f(interval(x));
  cout << "f=" << iy << ", width = " <<  width(iy) << endl;
    // Non-rigorous computations using multiprecision arithmetics 
  cout.precision(60);
  MpFloat::setDefaultPrecision(200);  // with 200 mantissa bits
  cout << "f=" << f(MpFloat(x)) << endl;
    // Rigorous computations using multiprecision interval arithmetics 
  MpInterval mpfx = f(MpInterval(x));
  cout << "f=" << mpfx << ",\nwidth = " <<  width(mpfx) << endl;
}
/* Output:
f=2.0811579830642382
f=[2.0811579830642328, 2.0811579830642457], width = 1.2878587085651816e-14
f=2.08115798306423807427767737012487687920820837024853353219391 
f=[2.08115798306423807427767737012487687920820837024853353219390 ,2.08115798306423807427767737012487687920820837024853353219392 ],
width = 7.46761833343337004857287686453614908870830260246540055972134e-60
*/
\end{lstlisting} 	

On top of these four basic arithmetic types the \capd library builds data structures such as vectors, matrices, hessians, jets (truncated Taylor series) and algorithms for manipulating them. Other data structures represent functions, solutions to ODEs or Poincar\'e maps, etc. Most of defined types use the following naming convention pattern
\begin{lstlisting}
[Prefix]ClassName
\end{lstlisting}
for example
\begin{lstlisting}
  DVector,   DMatrix,   DJet,   DMap,   DOdeSolver,   DPoincareMap, ...
 MpVector,  MpMatrix,  MpJet,  MpMap,  MpOdeSolver,  MpPoincareMap, ...
  IVector,   IMatrix,   IJet,   IMap,   IOdeSolver,   IPoincareMap, ...
MpIVector, MpIMatrix, MpIJet, MpIMap, MpIOdeSolver, MpIPoincareMap, ...
\end{lstlisting}

Prefixes \texttt{D} and \texttt{Mp} mean that the class is designed for non-rigorous computation based on \texttt{double} and \texttt{MpFloat} arithmetic types, respectively. Similarly, classes with prefixes \texttt{I} and \texttt{MpI} provide data structures and rigorous algorithms based respectively on \texttt{interval} and \texttt{MpInterval}. Whenever possible, we try to provide common interface for all kinds of data types and algorithms so that it is possible to switch between them if needed.

\subsection{Maps and their Taylor coefficients.}
%\lstinputlisting {code/vectalg.cpp}
One of the most important types is the class \texttt{[Prefix]Map} which represents a (possibly parameter dependent) map $$f_a:\mathbb R\times \mathbb R^n \ni (t,x_1,x_2,\dots,x_n) \to (f_1,f_2,\dots, f_m)\in \mathbb R^m,$$ where $a=(a_1, a_2,\dots, a_k)$ for some $k\geq 0$ is a vector of parameters. This class is usually used to define a generator of a discrete dynamical system or a vector field. The special time variable can be used to define non-autonomous vector fields.

This class provides also an easy to use interface for computation of higher order Taylor coefficients of the underlying map by means of automatic differentiation \cite{RallCorliss1996}. An instance of \texttt{Map} can be created by means of two constructors. If the map is given by a short formula it is convenient to parse the expression from a string with the following syntax
\begin{lstlisting}
IMap f("par:a1,a2,...,ak;time:t;var:x1,x2,...,xn;fun:f1,f2,...,fm;");
\end{lstlisting}
The sections \texttt{par} and \texttt{time} are optional. More complicated expressions can be defined as C++ functions with the signature
\begin{lstlisting}
void f(capd::autodiff::Node t,               // time variable
  capd::autodiff::Node in[], int dimIn,      // input variables x1,...,xn
  capd::autodiff::Node out[], int dimOut,    // output: function values
  capd::autodiff::Node params[], int noParam // parameters
);
\end{lstlisting}
and then sent to the constructor of \texttt{Map}. 
%Note that this function is used internally only to record computation pattern (directed acyclic graph) of an expression. In particular using relations (like comparison, inequalities) on input or any intermediate quantities does not make any sense and thus is not supported.
Below we present a short example illustrating the usage of both constructors and how the class can be used to compute values and derivatives of represented function.

\begin{lstlisting}
#include <iostream>
#include "capd/capdlib.h"
using namespace capd;
using namespace std;
/** Ikeda map is given by
 X = 1+u*(x*cos(r)-y*sin(r))
 Y =   u*(x*sin(r)+x*cos(r))
 where
 r = p - 6/(1+x^2+y^2) 
 and p,u are parameters */ 
void ikeda(capd::autodiff::Node /*t*/,       // unused time variable
  capd::autodiff::Node in[], int dimIn,      // input variables x,
  capd::autodiff::Node out[], int dimOut,    // output: function values X,Y 
  capd::autodiff::Node params[], int noParam // parameters: p,u
){
  capd::autodiff::Node r = params[0]-6./(1.+sqr(in[0])+sqr(in[1]));
  capd::autodiff::Node s = sin(r), c = cos(r);
  out[0] = 1+params[1]*(in[0]*c - in[1]*s);
  out[1] = params[1]*(in[0]*s + in[1]*c);  
}

int main(){ 
  IMap Henon("par:a,b;var:x,y;fun:y+1-a*x^2,b*x;");
  // Set parameter a=1.4. Note: it is not representable. 
  Henon.setParameter("a",interval(14)/10.);   
  // Set b\in [0.2,0.4], which contains the standard one b=0.3
  Henon.setParameter("b",interval(2, 4)/10.);
  IVector ix { {-1,2}, {0,1} };  // ix =[-1,2] x [0,1]
  IVector hx = Henon(ix);        // enclose Henon_{a,b}(ix)     
  // Enclose derivative D Henon_{a,b}(ix) 
  // for each point x in ix and for parameters a,b as above     
  IMatrix Dhx = Henon.derivative(ix);  

  int dimIn=2, dimOut=2, noParams=2, highestDerivative=5;
  IMap Ikeda(ikeda,dimIn,dimOut,noParams,highestDerivative);
  // Set parameters p=0.4, u=0.75, u is representable
  Ikeda.setParameters({interval(4)/10,interval(0.75)}); 
  IVector x {{1.,1.},{1.5,1.5}};  
  // Container for Taylor coefficient of 2-dimensional map up to order 4
  IJet jet(2,4);            
  Ikeda(ix, jet);   // Compute Taylor expansion of Ikeda map at point ix
  cout << jet.toString() << endl;    // and print in human readable form
  // Access to first order derivatives 
  cout << jet(0,0) << ", " << jet(0,1) << endl;
  cout << jet(1,0) << ", " << jet(1,1) << endl;
  // jet(i,j,c) gives access to normalized derivative (*@$\frac{\partial^2 f_i}{\partial x_j\partial x_c}$@*)
  cout << jet(0,0,1) << ", " << jet(1,0,0) << endl; // (*@$\frac{\partial^2 f_0}{\partial x_0\partial x_1}$@*) and (*@$\frac{1}{2!}\frac{\partial^2 f_1}{\partial x_0^2}$@*)
  // jet(i,j,c,k) gives access to normalized derivative (*@$\frac{\partial^3 f_i}{\partial x_j\partial x_c\partial x_k}$@*)
  cout << jet(1,0,0,1) << endl; // (*@$\frac{1}{2!}\frac{\partial^3 f_1}{\partial x_0^2\partial x_1}$@*)
  // For higher order Taylor coefficients use Multiindex notation
  cout << jet( Multiindex({2,2})) << endl; // Access to vector (*@$\frac{1}{2!2!}\frac{\partial^4 f}{\partial x_0^2\partial x_1^2}$@*)
}


\end{lstlisting}

\subsection{Solving initial value problems.}

Algorithms which solve initial value problems (IVPs) are split between three groups of classes.

\noindent\textbf{One-step solvers.} The first group consists of
\begin{lstlisting}
[Prefix]OdeSolver, [Prefix]CnOdeSolver
\end{lstlisting}
The above classes provide algorithms for one-step integration of ODEs. The class \texttt{OdeSolver} is optimized for $\mathcal C^0$ and $\mathcal C^1$ integration, that is solutions to IVPs and/or associated first order variational equations. The second class \texttt{CnOdeSolver} can integrate higher order variational equations as well. An example of its usage will be given in Section~\ref{sec:cr}.

\noindent\textbf{Long-time integration}. The above one-step methods are in general not recommended for direct usage. The second group of classes is built on top of \texttt{[Cn]OdeSolver}, that is
\begin{lstlisting}
[Prefix]TimeMap, [Prefix]CnTimeMap,
[Prefix]PoincareMap, [Prefix]CnPoincareMap.
\end{lstlisting}
Class \texttt{[Cn]TimeMap} combines a one-step solver with automatic step control strategies to compute trajectory segment over (usually) large time range. If integration time is not given explicitly but is determined by reaching certain Poincar\'e section, then one should use \texttt{[Cn]PoincareMap}. This class provides algorithms for computation of Poincar\'e maps and their derivatives. In the \capd library a Poincar\'e section is always defined as the set of zeroes of a smooth scalar-valued function $S:\mathbb R^m \to \mathbb{R}$ and realized by classes
\begin{lstlisting}
[Prefix]NonlinearSection, [Prefix]AffineSection, [Prefix]CoordinateSection.
\end{lstlisting}
The most general nonlinear case is covered by \texttt{NonlinearSection}. The library provides also computationally more efficient \texttt{AffineSection}, where the section is a hyperplane given by the normal vector $n\in \mathbb R^{m}$ and translation $c\in \mathbb R^m$, that is $S(x)=\left\langle n, x-c\right\rangle$. The last, and very often used type of section is \texttt{CoordinateSection}, where $S(x)=x_i - c$ for some $i\in1,\dots,m$ and a constant $c\in \mathbb R$.

\noindent\textbf{Sets and their propagation}. The third group consists of classes which specify different ways of representation of initial conditions and their propagation along trajectories. In rigorous computations special care should be taken on how intermediate results are represented. When a set of initial condition is
propagated by a dynamical system and on each step the image is bounded by an interval vector (product of intervals), then typically we observe the wrapping effect that leads to huge overestimation. On the other hand when the image is bounded by some non-linear shape, e.g given by multidimensional polynomials (like in the case of Taylor models \cite{BerzMakino1999}), then the result is more accurate but the computational cost increases rapidly with the dimension and degree of the polynomial.

In the \capd library the sets are represented (see \cite{MrozekZgliczynski2000}) as parallelepipeds, doubletons and tripletons. These strategies provide good compromise between speed and accuracy as shown in \cite{MIYAJI201634}. To choose appropriate set representation there are several factors to consider:
\begin{itemize}
    \item What set geometry will bring good compromise between speed and accuracy? The two that are usually the most efficient are
    \begin{itemize}
        \item \texttt{Rect2} - doubleton representation of the form $x + C*r0 + Q*q$ where
        $x,q,r0$ are interval vectors ($x$ is a point interval vector) and
        $C,Q$ are interval matrices, with $Q$ close to orthogonal.
        \item \texttt{Tripleton} - a subset of $ \mathbb R^m $ in the form $ x + C*r0 + \mathrm{intersection}(B*r,Q*q) $ where	
        $ x,q,r,r0$ are interval vectors ($ x $ is a point interval vector)
        and $ C,B,Q$ are interval matrices, with $ Q$ close to orthogonal.
    \end{itemize}
    \item What order of derivatives with respect to initial condition do we need? It is indicated by the prefix: \texttt{C0} sets enclose only the trajectory, \texttt{C1} sets enclose also first order derivatives with respect to initial conditions and \texttt{Cn} sets are used to store jets of flow up to given order, which has to be specified at the set construction.
    \item What numerical method should be used to propagate the set? There are two main groups of methods implemented in \capd: one based on the Taylor method and second based on the Hermite-Obreshkov (explicit-implicit) formula \cite{NedialkovJackson1998}. The infix \texttt{HO} indicates that the Hermite-Obreshkov method is requested.
    \item Is default double precision enough? If not add \texttt{Mp} prefix to compute with arbitrary precision (paying appropriate computational cost).
\end{itemize}
Summarizing, the names of data structures which represent initial conditions for ODEs follow the pattern
\begin{lstlisting}
[Mp]Cx[HO]GeometrySet
\end{lstlisting}
Not all combinations of components are implemented (please consult documentation for the full list of supported set representations), but the above pattern helps to encode a set representation type. For example
\begin{itemize}
\item \texttt{C0HOTripletonSet} stores $\mathcal C^0$ information only using tripleton representation and is propagated by the Hermite-Obershkov method \cite{NedialkovJackson1998},
\item \texttt{C1HORect2Set} stores points on the trajectory and first order derivatives with respect to initial conditions both in the form of doubletons and uses Hermite-Obreshkov method \cite{WalawskaWilczak2016} for their propagation.
\item \texttt{MpCnRect2Set} stores values and all derivatives up to given order in form of doubletons with \texttt{MpInterval} coefficients and propagates them by the Taylor method.
\end{itemize}

The following short code illustrates basic usage of the above three groups of classes.
\begin{lstlisting}
#include <iostream>
#include "capd/capdlib.h"
using namespace capd;
using namespace std;

int main(){
  cout.precision(11);
  // define vector field, solver and classes for long time integration
  IMap pendulum("par:w,pi;time:t;var:x,dx;fun:dx,cos(pi*t)-w*dx-sin(x);"); 
  pendulum.setParameters({interval(1)/100,interval::pi()});
  IOdeSolver solver(pendulum, 20);
  ITimeMap tm(solver);
  ICoordinateSection section(2,0); // two variables, index of x=0 is 0
  // alternatively, we could here use:
  // INonlinearSection section("var:x,y;fun:x;");
  // IAffineSection section(IVector(2),IVector({1.,0.})); //(origin,normal)
  IPoincareMap pm(solver,section);
  
  // integrate set [1,1]\times [0.25,0.5] over T=1
  C0HOTripletonSet s0( IVector({{1.,1.}, {0.25,0.5}}) );
  cout << "y=" << tm(1.,s0) << endl;

  // integrate point (1,0.5) and variational equations over T=100 and then  
  // print image and solution to variational equation (monodromy matrix)
  C1Rect2Set s1( IVector({{1.,1.}, {0.5,0.5}}) );
  cout << "u=" << tm(100.,s1) << endl;
  cout << "D=" << (IMatrix)s1 << endl;

  // continue integration of s1 until it reaches the section x=0
  IMatrix DP(2,2);
  IVector P = pm(s1,DP);
  cout << "P=" << P << endl;
  cout << "DP=" << DP << endl;
  // After computation P = pm(s1,DP) the set s1 is often far from section.
  cout << "s1=" << (IVector)s1 << endl;  
  cout << "Ds1=" << (IMatrix)s1 << endl;  
}
/* Output (rounded to 11 digits):
y={[1.0115905491, 1.2463432271],[-0.6386533516, -0.42846233913]}
u={[0.38323812527, 0.38323812528],[0.50996214875, 0.50996214877]}
D={
{[-4.3182170195, -4.3182170176],[-3.3700059498, -3.3700059482]},
{[4.3928041834, 4.3928041856],[3.3430225316, 3.3430225334]}
}
P={[-1.8354535396e-11, 1.8354535396e-11],[-0.34402409982, -0.34402409981]}
DP={
{[6.2681785123, 6.2681785151],[4.8172609505, 4.8172609528]},
{[-0.43855824706, -0.43855824587],[-0.27979305979, -0.27979305885]}
}
s1={[-0.063583452514, -0.063583452497],[-0.37811336848, -0.37811336848]}
Ds1={
{[6.0890412578, 6.0890412606],[4.6898120176, 4.6898120199]},
{[-1.5492512746, -1.5492512733],[-1.1344125491, -1.1344125481]}
} */

\end{lstlisting}

\subsection{The role of coordinate systems in integration of ODEs and computation of Poincar\'e maps}
The data structures, which represent initial conditions for ODEs provide constructors, that allow to set $x0,C,r0$ in both doubleton and tripleton representation. A proper usage of them can significantly improve obtained bounds. We will illustrate this issue with two suggestive examples. 

Consider the pendulum equation $x''=-\sin(x)$ and the following IVP: $x(0)\in[2,3]$ and $x'(0)=5-x(0)$, that is a line segment joining points $(2,3)$ and $(3,2)$ in the phase space. The following short code shows huge difference between two bounds on $(x(2),x'(2))$, depending on how this line segment is initially represented. The set \texttt{s1} represents initial conditions as 
$$
x0+C*r0 = \begin{bmatrix}2.5 \\2.5\end{bmatrix} + \begin{bmatrix}1 & 1 \\ -1 & 1\end{bmatrix}\begin{bmatrix}[-0.5,0.5]\\ [0,0]\end{bmatrix},
$$ 
while \texttt{s2} wraps it to the smallest interval vector, which contains this line segment, that is  $(x(0),x'(0))\in[2,3]^2$.

\begin{lstlisting}
#include <iostream>
#include "capd/capdlib.h"
using namespace capd;

int main(){
  IMap vf("var:x,dx;fun:dx,-sin(x);");
  IOdeSolver solver(vf,10);
  ITimeMap tm(solver);
  IVector x0({2.5,2.5}), r0({{-0.5,0.5},{0.,0.}});
  IMatrix C({{1.,1.},{-1.,1.}});
  C0HOTripletonSet s1(x0,C,r0), s2(x0+C*r0);
  std::cout << "s1(t=2): " << tm(2.,s1) << std::endl;
  std::cout << "s2(t=2): " << tm(2.,s2) << std::endl;
}
/** Output (rounded to 6 digits):
s1(t=2): {[7.821, 8.33752],[2.36664, 3.07598]}
s2(t=2): {[4.4177, 11.7408],[-0.930075, 6.3727]} */
\end{lstlisting}

Setting coordinate system on the Poincar\'e section also matters. Let us consider the Lorenz system \cite{Lorenz1963}
$$
x'= 10(y-x),\quad y'=x(28-z),\quad z'=xy-\frac{8}{3}z.
$$
Let us fix the Poincar\'e section $\Pi=\left\{(x,y,z) : z=27\right\}$ and denote by $\PM:\Pi\to\Pi$ the corresponding Poincar\'e map. We will use coordinates $(x,y)$ to describe points on $\Pi$. Define $\alpha = 7\pi/18$ and $$Q_\alpha=\begin{bmatrix}\cos\alpha & -\sin\alpha \\ \sin\alpha & \cos\alpha\end{bmatrix}.$$

A computer-assisted proof of the existence of chaotic dynamics in the Lorenz system given by Galias and Zgliczy\'nski \cite{GaliasZgliczynski98} required in particular, that the following inequality holds true: 
\begin{equation}\label{eq:lorenzInequality}
\left|\pi_y Q_\alpha^{-1}\PM^2(u)\right|<3.6 \quad \text{for}\quad  u=Q_\alpha\cdot(s,0),\ s\in[0.625,0.675].
\end{equation}
The following program illustrates the difference between multiplication by $Q_\alpha^{-1}$ after computation of Poincar\'e map, that is $Q_\alpha^{-1}\left(\PM^2(u)\right)$ and computation of $\left(Q_\alpha^{-1}\PM^2\right)(u)$ in a single routine. We see that the second estimate is much tighter and, in particular, the requested inequality (\ref{eq:lorenzInequality}) is validated.

\begin{lstlisting}
#include <iostream>
#include "capd/capdlib.h"
using namespace capd;

int main(){
  interval alpha = 7.*interval::pi()/18;
  interval c = cos(alpha), s = sin(alpha), z = 0.;
  IMatrix Q ({{c,-s,z},{s,c,z},{z,z,interval(1.)}});
  IMatrix invQ = transpose(Q);
  
  IMap lorenz("par:s;var:x,y,z;fun:10*(y-x),x*(28-z)-y,x*y-s*z;");
  lorenz.setParameter(0,interval(8)/3);
  IOdeSolver solver(lorenz,30);        // ODE integrator
  ICoordinateSection section(3,2,27.); // section is z=27
  IPoincareMap pm(solver,section);
  
  IVector r0(interval(625,675)/1000,-interval(36)/10,27);
  C0HOTripletonSet s1(IVector(3),Q,r0), s2 = s1;
  interval returnTime;
  // here we print (*@ $\pi_yQ^{-1}\left(\PM^2(x0+C*r0)\right)$  @*)
  std::cout << "y1: " << (invQ*pm(s1,2))[1] << std::endl;
  // here we print (*@ $\pi_y\left(Q^{-1}\PM^2\right)(x0+C*r0)$  @*)
  std::cout << "y2: " << pm(s2,IVector(3),invQ,returnTime,2)[1];
}
/** Output (rounded to 6 digits):
y1: [-4.16179, -0.129486]
y2: [-3.18066, -1.19211] */
\end{lstlisting}

\section{$\mathcal C^0$ solver and its applications}\label{sec:c0}
Topological methods are powerful and inexpensive in comparison to methods requiring estimates on derivatives. In this section we present two case studies:
\begin{itemize}
    \item the existence of symmetric periodic orbits in the Michelson system  and
    \item the existence of an attractor for the R\"ossler system.
\end{itemize}
\subsection{Periodic orbits in the Michelson system}
The Michelson system \cite{Mi} is a 3D system
\begin{equation}\label{eq:Michelson}
x'=y,\qquad y'=z,\qquad z'=c^2-y-x^2/2.
\end{equation}
reversible with respect to the involution
$$\mathcal R:(x,y,z) \to (-x,y,-z).$$
The above symmetry maps trajectories onto trajectories of the system but reverses the time, that is 
$\mathcal R(\phi(t,u)) = \phi(-t,\mathcal R(u))$
whenever $\phi(t,u)$ exists. A trajectory of (\ref{eq:Michelson}) is said to be $\mathcal R$-symmetric if it is invariant under this symmetry. 

\begin{figure}
    \centerline{\includegraphics[width=.95\textwidth]{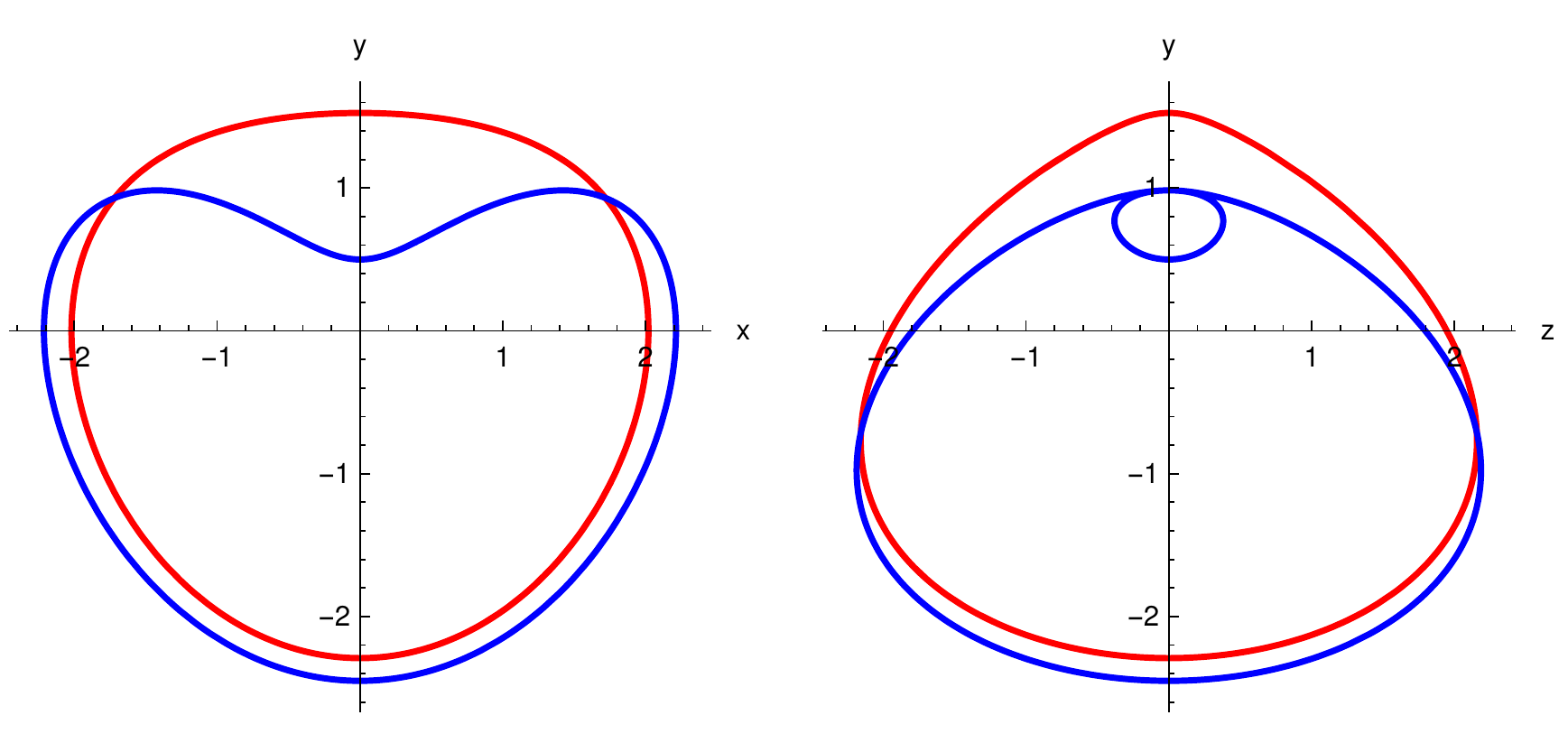}}
    \caption{Projections of two observed $\mathcal R$-symmetric periodic orbits for the system (\ref{eq:Michelson}) with the parameter $c=1$.\label{fig:twoper}}
\end{figure}

Let us define a Poincar\'e section $$\Pi = \left\{(0,y,z): y,z\in \mathbb R\right\}$$ and denote by
\begin{equation}\label{eq:MichelsonPM}
\PM_c\colon\Pi\to\Pi
\end{equation}
 a family of Poincar\'e maps parametrized by $c>0$. Note that we allow intersection of trajectories with the section $\Pi$ in both directions. It is easy to see that $\PM_c$ is reversible with respect to the involution $\mathcal R(y,z)=(y,-z)$, that is $\mathcal R\circ \PM = \PM^{-1}\circ \mathcal R$ --- see \cite{Wi1}. In \cite{TROY1989} an analytic proof of the existence of two $R$-symmetric periodic orbits for the system (\ref{eq:Michelson}) is given --- see Figure~\ref{fig:twoper}. Here we extend this result to a range of parameters.

\begin{theorem}
    For all parameter values $c\in C:=[1-1/128,1+1/128]$ the system (\ref{eq:Michelson}) has at least two $R$-symmetric periodic solutions.
\end{theorem}
\begin{proof}
    In order to prove the existence of a family of symmetric period-two points for $\PM_c$ in the set $\{0\}\times Y\times \{0\}$ parametrized by $c\in C$ it suffices to show that
    \begin{itemize}
        \item $\PM_c$ is defined and continuous on $\{0\}\times Y\times \{0\}$ for $c\in C$ and
        \item $\pi_z\PM_c(0,\min Y,0)$ and $\pi_z\PM_c(0,\max Y,0)$ have opposite signs for all $c\in C$.
    \end{itemize}
Then for each $c\in C$ there is a $y_c^*\in Y$ such that $\PM_c(0,y_c^*,0)=(0,\tilde y_c,0)$ for some $\tilde y_c\in\mathbb R$ and the result follows from the reversibility of $\PM_c$. The following program checks the above set of inequalities for two disjoint subintervals $Y_1, Y_2$ of positive semi-axis. Additionally, the program checks that
$\pi_y\PM_c(0\times y\times 0)<0$ for all $c\in C$ and $y\in Y_1\cup Y_2$, which implies
$$
    \left( \pi_y\PM_c(\{0\}\times Y_1\times \{0\})\cup \pi_y\PM_c(\{0\}\times Y_2\times \{0\})\right)\cap (Y_1\cup Y_2) = \emptyset
$$
and thus the two families of periodic points in $Y_1$ and $Y_2$ are different. The program executes within less than 1 second on a laptop-type computer.
\end{proof}
\begin{lstlisting}
/* Proof of symmetric periodic orbits in the Michelson system. */ 
#include <iostream>
#include "capd/capdlib.h"
using namespace capd;
/**
 * @param c - parameter range
 * @param Y - interval on the y-axis
 * @param n - numer of grid elements in Y
 */ 
void validateSymPO(interval c, interval Y, int n){
  IMap vf("par:c;var:x,y,z;fun:y,z,c^2-y-0.5*x^2;");
  IOdeSolver solver(vf, 15);        // order of the solver    
  ICoordinateSection section(3, 0); // Poincare section x=0
  IPoincareMap pm(solver, section);
  vf.setParameter("c",c); 
  // We split Y into n subintervals and check that
  // the Poincare map is defined on (0,Y ,0). What the result is
  // is not crucial, as long as the image is in the y<0 halfplane.
  interval g = (Y.right()-Y.left())/n;
  for(int i=0;i<n;++i){
    C0HOTripletonSet s({0.,Y.left() + interval(i,i+1)*g,0.});
    assert( pm(s)[1] < 0.);  
  }
  // Compute P(0,min(Y),0) and P(0,max(Y),0)
  C0HOTripletonSet s1(IVector({0.,Y.left(),0.}));
  C0HOTripletonSet s2(IVector({0.,Y.right(),0.}));
  IVector r1 = pm(s1), r2 = pm(s2);
  // Check that z-coordinate changes the sign
  std::cout << "validated? " << ( r1[2]*r2[2]<0 ) << ": " 
            << r1[2] << ", " << r2[2] << std::endl;
}
int main(){
  std::cout.precision(4);
  interval I(-1,1);
  // Call the algorithm with two approximate periodic points. 
  // Parameter c=1, validate orbits with accuracy 1e-13.
  validateSymPO(1., 1.5259617305036892  + 1e-13*I, 1);
  validateSymPO(1., 0.50002564853520548 + 1e-13*I, 1);
  // Repeat for parameter range c \in [1-1/128,1+1/128].
  validateSymPO(1+I/128, 1.5259617305036892  + 2e-1*I, 2);
  validateSymPO(1+I/128, 0.50002564853520548 + 2e-1*I, 4);  
  return 0;
}
/* Ouptut (rounded to 4 digits, change of the signs is relevant):
validated? 1: [1.414e-13, 2.466e-13], [-2.48e-13, -1.444e-13]
validated? 1: [-4.212e-13, -2.477e-13], [2.331e-13, 4.036e-13]
validated? 1: [0.1046, 0.6294], [-0.6001, -0.1967]
validated? 1: [-1.314, -0.09417], [0.03929, 1.153] */
\end{lstlisting}

\subsection{Attractor in the R\"ossler system}\label{sec:RosslerAttractor}
\begin{figure}[htbp]
    \centerline{
        \includegraphics[height=2.2in]{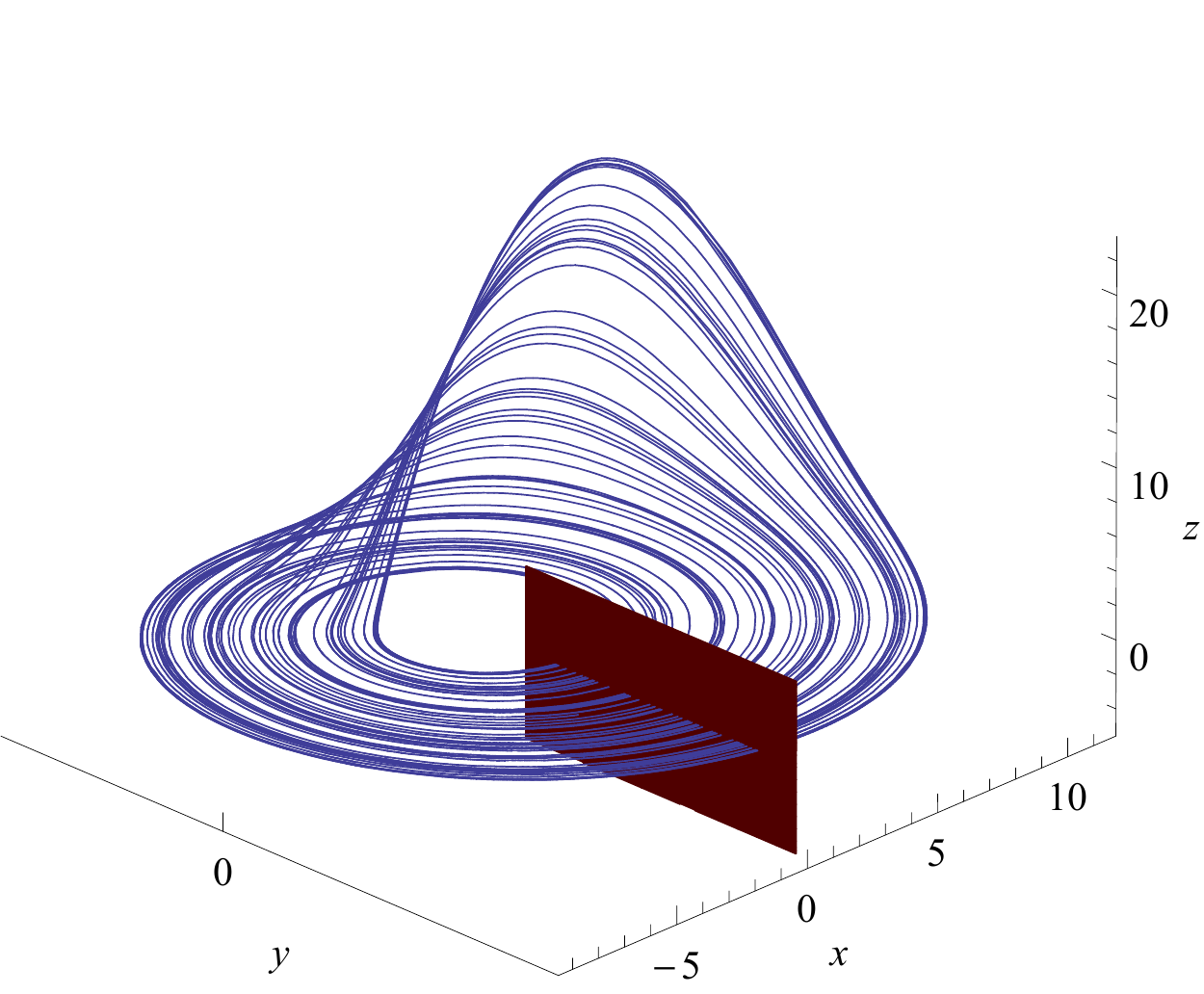}\quad
        \includegraphics[height=2.2in]{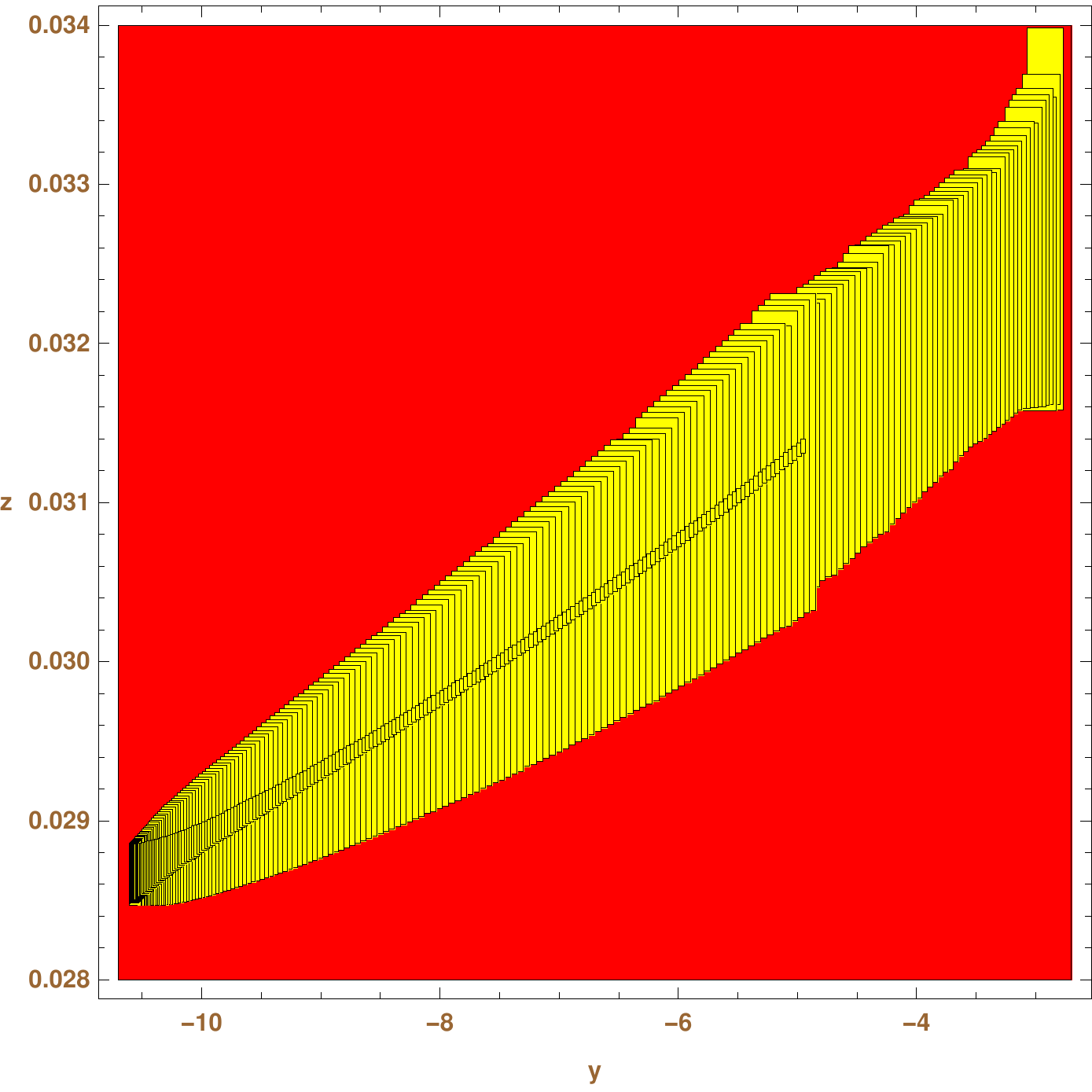}
    }
\caption{(Left) Observed chaotic attractor of the system (\ref{eq:Rossler}) with $b=0.2$ and $a=5.7$. (Right) Plot of trapping region $W$ and rigorous enclosure of $\PM_{a,b}(W)$ for the system with parameters $b=0.2$ and $a=5.7$.\label{fig:Rossler}}
\end{figure}

Simulation shows that for a wide range of parameter values the R\"ossler system \cite{Rossler1976}
\begin{equation}\label{eq:Rossler}
x'=-(y+z),\qquad y'=x+by,\qquad z'= b+z(x-a)
\end{equation}
possesses a chaotic attractor --- see Fig.~\ref{fig:Rossler}. To the best of our knowledge, the first proof that for classical parameter values there is a trapping region for the attractor was given in \cite{WalawskaWilczak2019}.
To present the details of the computer assisted proof, let us define a Poincar\'e section and the  corresponding Poincar\'e map by
\begin{equation}\label{eq:RosslerPM}
\begin{array}{lll}
\Pi &=& \left\{(0,y,z): y,z\in \mathbb R,x'>0\right\}, \\
\PM_{a,b}&:&\Pi\to\Pi.
\end{array}
\end{equation}
\begin{theorem}[\cite{WalawskaWilczak2016}]\label{thm:RosslerTrappingRegion}
    For $a=5.7$ and $b=0.2$ the set
    \begin{equation}\label{eq:RosslerTrappingRegion}
    W = Y\times Z := [-10.7,-2.7] \times [0.028,0.034]
    \end{equation}
    is forward invariant for $\PM_{a,b}$. In particular, it contains compact, connected invariant set $\mathcal A=\bigcap_{n>0}\PM_{a,b}^n(W)$.
\end{theorem}
\begin{proof}
    First observe, that $x'=-(y+z)\geq 2.7-0.034>0$ for all $(y,z)\in W$ and hence $W\subset \Pi$.
    Inclusion $\PM_{a,b}(W)\subset W$ can be checked in direct computation. The set $W$ is subdivided uniformly (for simplicity of the program) $W\subset \bigcup_{i=1}^{N} Y_i\times Z$ and then inclusion $\PM_{a,b}(Y_i\times Z)\subset W$ is checked for all $i=1,\ldots,N$. The program executes within less than 1 second on a laptop-type computer.
\end{proof}
\begin{lstlisting}
#include <iostream>
#include "capd/capdlib.h"
using namespace capd;
int main(){
  IMap vf("par:a,b;var:x,y,z;fun:-(y+z),x+b*y,b+z*(x-a);");
  vf.setParameter("a",interval(57)/10);
  vf.setParameter("b",interval(2)/10);
  IOdeSolver solver(vf, 20);
  ICoordinateSection section(3, 0);   // section is given by x = 0 and
  IPoincareMap pm(solver, section, poincare::MinusPlus);  // x'> 0
  // Coordinates of the trapping region W = Y\times Z
  interval Y = interval(-107,-27)/10; // Y=[-10.7,-2.7] 
  interval Z = interval(28,34)/1000;  // Z=[0.028,0.034]
  // Subdivide uniformly Y onto N subintervals
  const int N = 150;
  bool result = true;
  for (int i = 0; i < N and result; ++i) {
    IVector Wi ({0., Y.left() + diam(Y)*interval(i,i+1)/N, Z});
    C0HOTripletonSet s(Wi);
    IVector u = pm(s);
    result = result and subset(u[1],Y) and subset (u[2],Z);
  }
  std::cout << "Trapping region validated? : " << result << std::endl;
  return 0;
}
\end{lstlisting}

In the next section we will show that the attractor $\mathcal A=\bigcap_{n>0}\PM_{a,b}^n(W)$ contains a chaotic and uniformly hyperbolic invariant set.

\section{$\mathcal C^1$ solver and its applications}\label{sec:c1}
By a rigorous $\mathcal C^1$-algorithm we mean an algorithm capable of computing bounds for the following system of ODEs
\begin{equation*}
x'(t)= f(t,x(t)),\qquad V'(t) = D_xf(t,x(t))\cdot V(t)
\end{equation*}
for some initial conditions $x(0)\in[x_0]\subset \mathbb R^n$ and $V(0)\in [V_0]\subset \mathbb R^{n\times n}$. In principle, any $\mathcal C^0$ solver is capable of doing this task. Taking into account special structure of this system of equations, one can design an algorithm of complexity $O(n^3)$  which is much faster than the direct application of the $\mathcal C^0$ solver which has complexity $O(n^6)$. Such an algorithm was proposed in \cite{Zgliczynski2002} and later improved in \cite{WalawskaWilczak2016}. Both versions are very powerful tools for studying hyperbolic-like properties of dynamical systems, such as verification of periodic orbits and their stability \cite{KapelaZgliczynski2003,KapelaSimo2007,BarrioRodriguezBlessa2012}, connecting orbits \cite{WilczakZgliczynski2003,Wilczak2009abundance} and hyperbolic attractors \cite{Wilczak2010}. Here we present a few short, yet non-trivial examples:
\begin{itemize}
    \item verification of the existence of a solution to a boundary value problem,
    \item verification of the existence of hyperbolic periodic solutions with very high localization accuracy,
    \item verification of the existence of a  hyperbolic chaotic set.
\end{itemize}

\subsection{Boundary value problem.}
In this section we show that $\mathcal C^1$ algorithms can be used to solve boundary value problems. As an example, we reproduce the result by Nakao \cite{NAKAO1992489}.
\begin{theorem}[\cite{NAKAO1992489}]
    The equation
    \begin{equation}\label{eq:Nakao}x''=-0.1x-0.1x^3-0.4464\cos t\end{equation}
    has a solution satisfying $x'(0)=x'(2\pi)=0$.
\end{theorem}
\begin{proof}
    The proof in \cite{NAKAO1992489} is also computer-assisted but relies on solving a zero-finding problem in some infinite-dimensional functional space. Here we propose a direct approach, as we have tools capable to compute derivatives of ODEs with respect to initial conditions. Denote by $\varphi(t,t_0,x_0,x'_0)=(\varphi_x(t,t_0,x_0,x_0'),\varphi_{\dot x}(t,t_0,x_0,x_0'))$ a solution of the initial value problem $x(t_0)=x_0, x'(t_0)=x_0'$ for (\ref{eq:Nakao}). It is easy to see that the zeroes of
    $$F(x) = \varphi_{\dot x}(2\pi,0,x,0) = 0$$
    correspond to the solutions of the boundary value problem we are looking for. The following program checks, by means of the interval Newton operator (\ref{eq:NewtonOperator}), that the function $F$ has a zero at some $x_*$ with $|x_*+0.5072|\leq 10^{-4}$. The program executes within less than 1 second on a laptop-type computer.
\end{proof}
\begin{lstlisting}
/** Example of solving BVP: x'(0)=x'(2pi)=0 **/
#include <iostream>
#include "capd/capdlib.h"
using namespace capd;
using namespace std;
int main(){
  IMap f("par:a;time:t;var:x,dx;fun:dx,-x*(1+x^2)/10 + a*cos(t);");
  f.setParameter("a",interval(4464)/10000);
  IOdeSolver solver(f,20);  // ODE integrator
  ITimeMap tm(solver);      // class for long time integration

  IVector u0({-0.5072,0.});
  IVector r({interval(-1e-4,1e-4),0.});
  C0HOTripletonSet s0(u0);   
  C1HORect2Set s(u0,r);  // represent set s = u0 + Id*r
  // integrate set and variational equation until T=2*pi 
  IVector y = tm(2.*interval::pi(),s0); // C0 computation  
  tm(2.*interval::pi(),s);              // C1 computation
  // solve equation F(r1) := proj_{x'}(phi(2pi,u0+(r1,0))) = 0
  interval N = - y[1]/((IMatrix)s)[1][0];
  cout << "(N,r1)=(" << N << "," << r[0] << ")"<< endl;
  cout << "subset(N,r)? = " << boolalpha << subset(N,r[0]) << endl;
  return 0;
}
/* Output (rounded to six digits):
(N,r1)=([-2.68037e-05, -2.67903e-05],[-0.0001, 0.0001])
subset(N,r)? = true
*/ 
\end{lstlisting}

\subsection{High localization accuracy bounds for periodic orbits in the R\"ossler system.}
In Section~\ref{sec:RosslerAttractor} we have proved that  system (\ref{eq:Rossler}) has an attractor. Here we will prove the existence of three periodic orbits on this attractor with very high localization accuracy.
\begin{theorem}
    Put
        \begin{eqnarray*}
        u_1 &=& (-8.3809417428298762873487630431,0.029590060630667102951494027735),\\
        u_2 &=& (-5.4240738226652043515673025463,0.031081210807876445187367377796),\\
        u_3 &=& (-6.2331586285379749515076479411,0.030640111658160569478006226700).
    \end{eqnarray*}
    There exist three hyperbolic periodic points $u_m^*=u_m+r_m$, $m=1,2,3$ for the Poincar\'e map (\ref{eq:RosslerPM}) of period $1,2,3$, respectively satisfying $\|r_m\|_1\leq 10^{-28}$ and the coordinates of $r_m$ are known with accuracy $10^{-54}$.
\end{theorem}
\begin{proof}
    The proof relies on properties of the interval Newton operator \cite{Neumaier1990}
    \begin{equation}\label{eq:NewtonOperator}
    N(f,x_0,X) = x_0 - [Df(X)]_I^{-1}f(x_0),
    \end{equation}
where by $[A]_I$ we mean an interval hull of the matrix $A$. It is well known \cite{Neumaier1990}, that if $X$ is a convex set, $x_0\in X$ and $N(f,x_0,X)\subset  X$ then the mapping $f$ has a unique zero in the set $X$. Moreover, this zero belongs to $N(f,x_0,X)$. The following program validates the existence of three periodic solutions by means of the interval Newton operator applied to the  function $f=\PM_{a,b}^m-\mathrm{Id}$ for $m=1,2,3$, depending on the orbit. It also prints the largest diameter of the components of $N(\PM_{a,b}^m-\mathrm{Id},u_m,u_m+R_m)-u_m$, where $R_m=10^{-28}\cdot [-1,1]^2$, which in each case $m=1,2,3$ is less than $10^{-54}$. This shows, that coordinates of $r_m\in N(\PM_{a,b}^m-\mathrm{Id},u_m,u_m+R_m)-u_m$ are known with accuracy $10^{-54}$. We would like to emphasize, that it is very easy to obtain much higher localization accuracy by either providing more digits for initial points or by iterating the interval Newton operator. Finally, the program computes bounds on the eigenvalues of the Poincar\'e map at periodic points which proves that they are all of saddle type. From these bounds it is also clear, that the three orbits are different.

The program uses high-precision version of the $\mathcal C^1$ ODE solver from the \capd library to obtain tiny bounds on $\PM_{a,b}$ and its derivative. The program executes within 20 seconds on a laptop-type computer.
\end{proof}
\begin{lstlisting}
#include <iostream>
#include "capd/mpcapdlib.h"
using namespace capd; 
using namespace std;
/**
 * @param (y,z) - approximate periodic point
 * @param e - radius of the set centred at (y,z)
 * @param n - period
 */ 
void po(MpFloat y, MpFloat z, int n, double e=1e-28){
 MpIMap vf("par:a,b;var:x,y,z;fun:-(y+z),x+b*y,b+z*(x-a);");
 vf.setParameter("a",MpInterval(57)/10);
 vf.setParameter("b",MpInterval(2)/10);
 MpIOdeSolver solver(vf,80);           // ODE integrator of order 80
 MpICoordinateSection section(3,0.);   // the section is x=0
 MpIPoincareMap pm(solver,section, poincare::MinusPlus);
 // Approximate periodic point and a ball around it.
 MpIVector u0({MpInterval(0.),y,z});
 MpIVector r({MpInterval(0.),MpInterval(-e,e),MpInterval(-e,e)});
 MpC0TripletonSet s0(u0);
 // Compute P^n(u0)-u0 and project it onto (y,z)
 MpIVector fu0( 2, (pm(s0,n) - u0).begin() + 1 ); 
 MpC1Rect2Set s1(u0,r); // represent s1 =  u0 + Id*r   
 MpIMatrix Dphi(3,3);
 // compute DP^n(u0+r) 
 MpIVector u = pm(s1,Dphi,n);
 MpIMatrix DP = pm.computeDP(u,Dphi);
 // projection of DP^n(u0+r)-Id onto 2D subspace
 MpIMatrix M({{DP[1][1]-1.,DP[1][2]},{DP[2][1],DP[2][2]-1.}});
 // enclose -(DP^n(u0+r)-Id)^{-1}*(P^n(u0)-u0)
 MpIVector N = - matrixAlgorithms::gauss(M,fu0);
 cout << boolalpha << "(validated?, accuracy) = (" 
      << subset(N,MpIVector(2,r.begin()+1)) << ", " << maxWidth(N) << ")";
 // compute bound on eigenvalues using an explicit formula
 MpInterval t = sqrt(4*DP[2][1]*DP[1][2] + sqr(DP[1][1]-DP[2][2])); 
 cout << "\neigenvalues=(" << 0.5*(DP[1][1]+DP[2][2]-t) 
      << "," << 0.5*(DP[1][1]+DP[2][2]+t) << ")" << endl;
}
int main(){
 MpFloat::setDefaultPrecision(200); // 200 manitissa bits
 po("-8.3809417428298762873487630431",".029590060630667102951494027735",1);
 po("-5.4240738226652043515673025463",".031081210807876445187367377796",2);
 po("-6.2331586285379749515076479411",".030640111658160569478006226700",3);
}
/* Output (rounded to 6 digits) of the program:
(validated?, accuracy) = (true, 8.10612e-56 )
eigenvalues=([-2.40396 ,-2.40395 ],[-1.28211e-14 ,-1.28210e-14 ])
(validated?, accuracy) = (true, 2.03877e-55 )
eigenvalues=([-3.51201 ,-3.51200 ],[-1.24264e-26 ,1.20278e-26 ])
(validated?, accuracy) = (true, 6.49223e-55 )
eigenvalues=([-2.34193 ,-2.34192 ],[-2.73947e-26 ,2.73947e-26 ])*/
\end{lstlisting}

\subsection{Uniformly hyperbolic chaotic invariant set}
In the last example of this section we would like to recall the result from \cite{WalawskaWilczak2016} about the existence of a uniformly hyperbolic and chaotic invariant set in the R\"ossler system (\ref{eq:Rossler}). From Theorem~\ref{thm:RosslerTrappingRegion} we know that the set $W=Y\times Z$  defined by (\ref{eq:RosslerTrappingRegion}) is a trapping region for the Poincar\'e map (\ref{eq:RosslerPM}) of the R\"ossler system (\ref{eq:Rossler}) for parameters values $a=5.7, b=0.2$. This implies the existence of a connected, compact invariant set $\mathcal A=\bigcap_{n>0}\PM_{a,b}^n(W)$. Here we present a proof that this attractor is non-trivial.
\begin{theorem}[\cite{WalawskaWilczak2016}]
    Let $l_M=-8.4$, $r_M = -7.6$, $l_N=-5.7$, $r_N=-4.6$ and define two subsets of $W$,
    \begin{equation*}
        M = [l_M,r_M]\times Z \quad\textrm{ and }\quad N = [l_N,r_N]\times Z.
    \end{equation*}
    Fix $a=5.7$, $b=0.2$ and denote $\PM=\PM_{a,b}$.
    Then the maximal invariant set for $\PM^2$ in $N\cup M$, denoted by $\mathcal H
        =\mathrm{inv}(\PM^2,N\cup M)\subset \mathcal A$, is uniformly hyperbolic; in
        particular it is robust under perturbations of the system. The dynamics of
        $\PM^2$ on $\mathcal H$ is chaotic in the sense that $\PM^2|_{\mathcal H}$ is
        conjugated to the Bernoulli shift on two symbols.
\end{theorem}
\begin{proof}
The proof relies on some partial results from~\cite{Zgliczynski1997,KokubuWilczakZgliczynski2007,Wilczak2010}. Semiconjugacy of $\PM^2|_{\mathcal H}$ to the Bernoulli shift is proved by means of the method of covering relations \cite{Zgliczynski1997}. It is sufficient to check the following geometric conditions
\begin{equation}\label{eq:covrel}
\begin{array}{lclc}
\pi_y\PM_{a,b}^2(y,z)&<& l_M \quad \text{for }\ (y,z)\in\{l_M\}\times Z,\\
\pi_y\PM^2(y,z)&>& r_N \quad \text{for }\ (y,z)\in\{r_M\}\times Z,\\
\pi_y\PM^2(y,z)&<& l_M \quad \text{for }\ (y,z)\in\{r_N\}\times Z, \\
\pi_y\PM^2(y,z)&>& r_N \quad \text{for }\ (y,z)\in\{l_N\}\times Z,
\end{array}
\end{equation}
where $\pi_y$ denotes the projection onto $y$ coordinate. Rigorous bounds on $\PM^2(\{l_M\}\times Z)$, $\PM^2(\{r_M\}\times Z)$, $\PM^2(\{l_N\}\times Z)$ and $\PM^2(\{r_N\}\times Z)$, returned by our
routine, are shown in Fig.~\ref{fig:chaos}.

\begin{figure}[htbp]
    \centerline{\includegraphics[width=.7\textwidth]{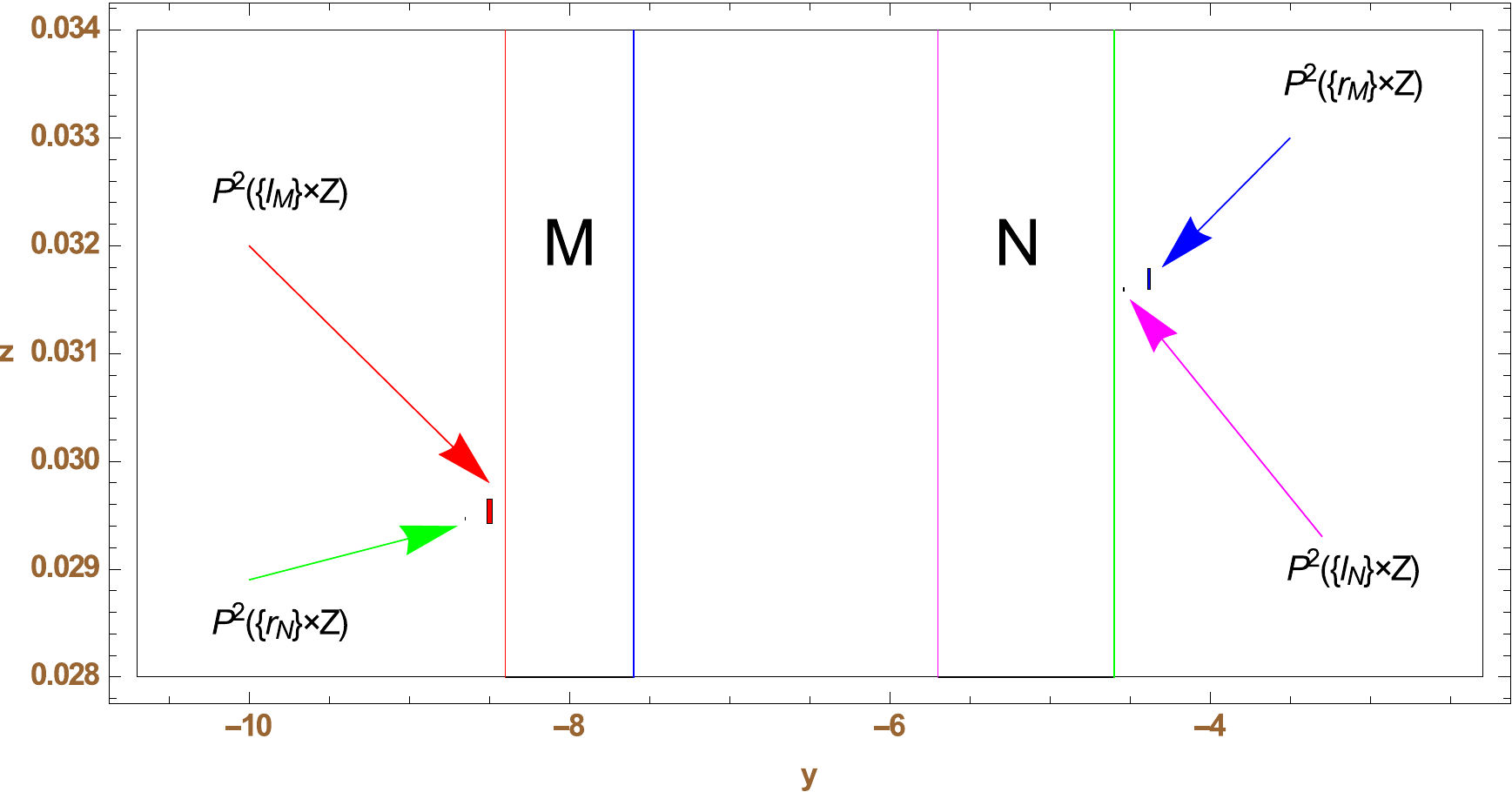}}
    \caption{The sets $M$ and $N$ and rigorous enclosures of the images of their
        exit edges --- see (\ref{eq:covrel}).\label{fig:chaos}}
\end{figure}
Hyperbolicity of $\mathcal  H$ is proved by means of the cone condition introduced in
\cite{KokubuWilczakZgliczynski2007}. Let $Q$ be a diagonal matrix $Q=\mathrm{Diag}(\lambda,\mu)$ with arbitrary coefficients satisfying $\lambda>0$ and $\mu<0$. It was shown in \cite{Wilczak2010} that if
for all $(y,z)\in N\cup M$ the matrix
\begin{equation}\label{eq:coneCondition}
DP^2(y,z)^T\cdot Q \cdot DP^2(y,z) -Q
\end{equation}
is positive definite, then the maximal invariant set for $P^2$ in $N\cup M$
is uniformly hyperbolic. The following program checks all necessary conditions with constants $\lambda=1$ and
$\mu=-100$. Some subdivision of sets was necessary to obtain sharp bounds on the derivative of $\PM^2$. The program executes within less than 2 seconds on a laptop-type computer.
\end{proof}
\begin{lstlisting}
#include <iostream>
#include "capd/capdlib.h"
using namespace capd;
using namespace std;

// z-coordinate of the trapping region
interval Z = interval(28,34)/1000;  // Z=[0.028,0.034]
// y-coordinates of sets M and N
interval My=interval(-84,-76)/10, Ny=interval(-57,-46)/10;

/// This routine checks cone-condition on the set Y\times Z
/// The interval Y is split into g pieces
bool checkCC(IPoincareMap& pm, interval Y, int g) {
  bool res = true;
  interval p = (Y.right()-Y.left())/g;
  IMatrix Dphi(3,3);
  // define quadratic form Q
  IMatrix Q({{0.,0.,0.},{0.,1,0.},{0.,0.,-100}});
  for (int i = 0; i < g and res; ++i) {
    // compute derivative of P^2 on a grid element
    C1Rect2Set s({0.,Y.left()+interval(i,i+1)*p,Z});
    interval returnTime; // not used, by required by syntax below 
    IVector y = pm(s, Dphi, returnTime, 2);
    IMatrix DP = pm.computeDP(y,Dphi);
    // check positive definitness by Sylvester criterion
    DP = transpose(DP)*Q*DP - Q;
    res = res and DP[1][1]>0 and (DP[1][1]*DP[2][2]-sqr(DP[1][2]))>0;
  }
  return res;
}
int main(){
  IMap vf("par:a,b;var:x,y,z;fun:-(y+z),x+b*y,b+z*(x-a);");
  vf.setParameter("a",interval(57)/10);
  vf.setParameter("b",interval(2)/10);
  IOdeSolver solver(vf, 20);
  ICoordinateSection section(3, 0); // section x=0, x'>0
  IPoincareMap pm(solver, section, poincare::MinusPlus);

  // Inequalities for the covering relations -- see (*@(\ref{eq:covrel})@*).
  cout << boolalpha;
  C0HOTripletonSet LM( IVector({0.,My.left(), Z}) ); 
  C0HOTripletonSet RM( IVector({0.,My.right(),Z}) ); 
  C0HOTripletonSet LN( IVector({0.,Ny.left(), Z}) ); 
  C0HOTripletonSet RN( IVector({0.,Ny.right(),Z}) ); 
  cout << "pi_y P^2(LM) < lM? " << ( pm(LM,2)[1] < My ) << endl;
  cout << "pi_y P^2(RM) > rN? " << ( pm(RM,2)[1] > Ny ) << endl;
  cout << "pi_y P^2(LN) > rN? " << ( pm(LN,2)[1] > Ny ) << endl;
  cout << "pi_y P^2(RN) < lM? " << ( pm(RN,2)[1] < My ) << endl;  
  // Check cone conditions.
  cout << "Cone condition on M? " << checkCC(pm,My,80) << endl;
  cout << "Cone condition on N? " << checkCC(pm,Ny,40) << endl;
}
\end{lstlisting}
\section{$\mathcal C^r$ solver and its application}\label{sec:cr}
\begin{figure}
    \centerline{\includegraphics[width=.7\textwidth]{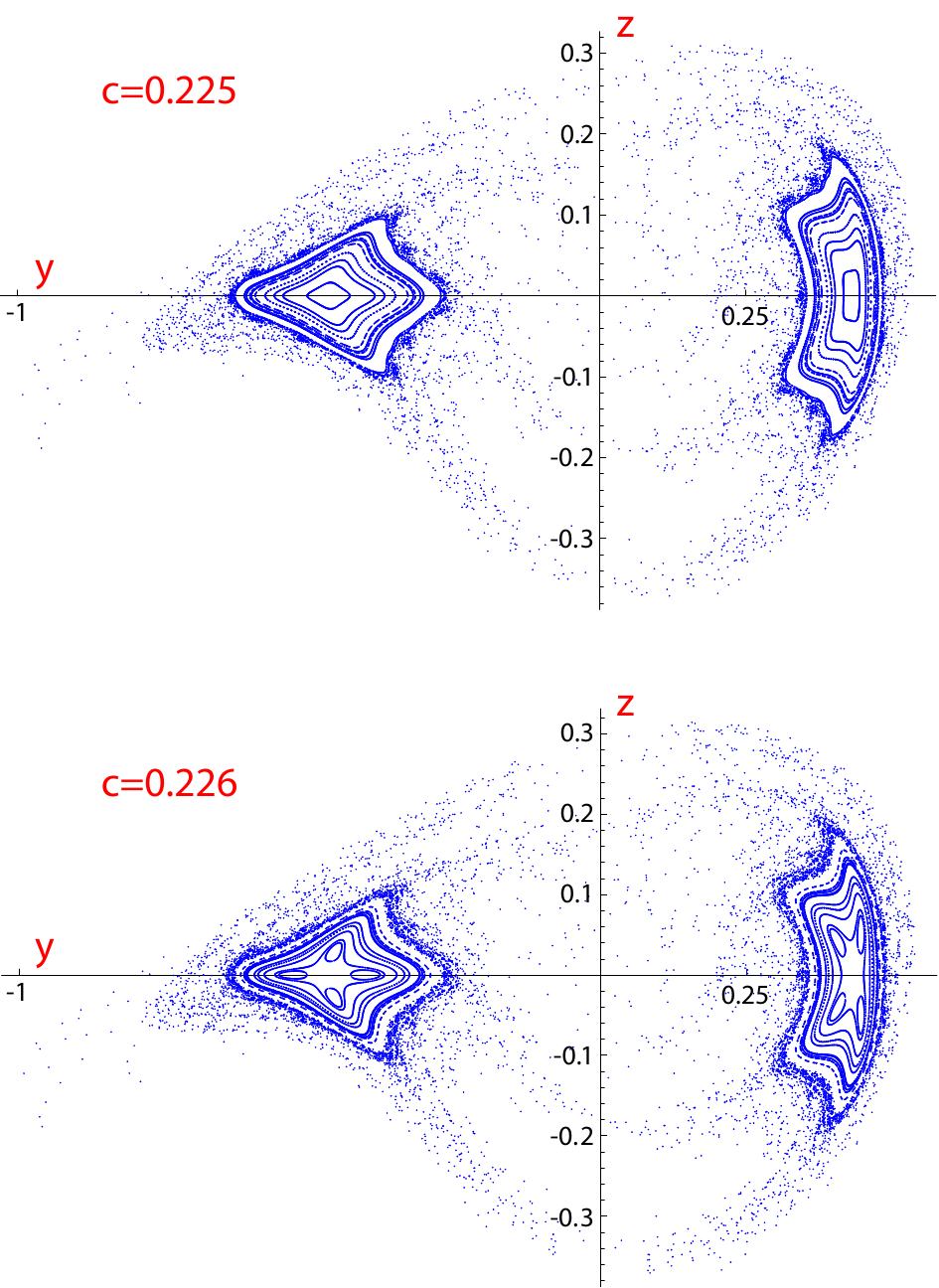}}
    \caption{Phase portrait of the Poincar\'e map (\ref{eq:MichelsonPM}) of the Michelson system (\ref{eq:Michelson}). It is an evidence that a family of elliptic periodic orbits crosses $1:4$ resonance when the parameter varies in the interval $c\in[0.225,0.226]$. The existence of such period quadrupling bifurcation was proved in \cite{WalawskaWilczak2019}.\label{fig:ellipticOrbit}}
\end{figure}

The \capd library offers a rigorous solver for higher order variational equations, that is
\begin{eqnarray*}
\frac{d}{dt}\phi(t,x) &=& f(t,x(t)),\\
\frac{d}{dt}D_x\phi(t,x) &=& D_xf(t,x(t))\cdot D_x\phi(t,x),\\
\frac{d}{dt}D_a\phi(t,x) &=& D_xf(t,x(t))D_a\phi(t,x) + h.o.t,
\end{eqnarray*}
where $D_a$ is the partial derivative operator with respect to a multiindex $a$ and $h.o.t.$ stands for higher order terms not written explicitly. Higher order derivatives are extremely useful in studying global and local bifurcations \cite{KokubuWilczakZgliczynski2007,WilczakZgliczynski2009focm,WalawskaWilczak2019,WilczakZgliczynski2009siads} as well as non-linear stability of elliptic periodic solutions \cite{KapelaSimo2017,WilczakBarrio2017}. Here we present one short example of application of $\mathcal C^r$  solver to study KAM tori near an elliptic periodic orbit in the Michelson system (\ref{eq:Michelson}). The existence of a wide branch of such orbits parametrized by $c$ was proved in \cite{WilczakBarrio2017}. Here we give a proof for just one parameter value close to $1:4$ resonance -- see Fig.~\ref{fig:ellipticOrbit}.
The following theorem is a special case of the result from \cite{WilczakBarrio2017}.
\begin{theorem}
    For the parameter value $c=0.226$ the Poincar\'e map (\ref{eq:MichelsonPM}) has a symmetric period-two point $u_*=(0,y_*,0)$, $|y_*-0.43407644067709|\leq 2\cdot 10^{-14}$, which is stable. That is, any neighbourhood $U\subset\mathbb R^3$ of the periodic trajectory $\mathcal O(u_*)$ contains a $2D$ invariant torus surrounding the orbit $\mathcal O(u_*)$ and separating the phase space.
\end{theorem}
\begin{proof}
We will use computational tools to validate the hypotheses of the classical result by Siegel and Moser \cite{SiegelMoser1971}. First, we have to check that the periodic orbit indeed exists. The interval Newton operator (\ref{eq:NewtonOperator}) applied to the scalar equation
    $$
    F(y) = \pi_z(\PM_c(0,y,0)) = 0
    $$
    validates the existence of a zero $y_*$ of $F$ satisfying $|y_*-0.43407644067709|\leq 2\cdot 10^{-14}$.

         Then we compute truncated Birkhoff normal form at the periodic orbit and we check the twist condition. If a certain coefficient in the normal form does not vanish then the existence of KAM tori and stability of periodic orbit follow from the theorem by Siegel and Moser. The following program executes within less than 1 second on a laptop-type computer.
\end{proof}
\begin{lstlisting}
/** Existence of invariant curves around an elliptic PO **/
#include <iostream>
#include "capd/capdlib.h"
#include "capd/normalForms/planarMaps.hpp"
using namespace capd;
using namespace std;
int main(){
  cout.precision(17);
  IMap vf("par:c;var:x,y,z;fun:y,z,c^2-y-0.5*x^2;",3);
  vf.setParameter("c",interval(226)/1000);
  ICnOdeSolver solver(vf,20);
  ICoordinateSection section(3,0); // x=0
  ICnPoincareMap pm(solver,section);           
  IMatrix DP(3,3);

  // Validate existence of a periodic point by the Newton method
  IVector u0({0,interval(43407644067709L)/1e+14,0.});
  IVector r = 1e-12*interval(-1,1)*IVector{0.,1,0.};
  C0TripletonSet s0(u0);
  IVector y0 =  pm(s0);
  C1Rect2Set s(u0,r); // represent s = u0 + Id*r
  IVector y = pm(s,DP);
  DP = pm.computeDP(y,DP);
  interval N = - y0[2]/DP[2][1];
  cout << "subset(N,r)? = " << boolalpha << subset(N,r[1]) 
       << ", N = " << N << endl;
  // Integrate 3rd order variational equations over the full period 2
  CnRect2Set S(u0+N,3);
  IJet jet(3,3,3);
  y = pm(S,jet,2);
  jet = pm.computeDP(jet); 
  // project onto 2-dim section
  IJet P(2,2,3);
  for(int j=0;j<=3;++j)
    for(int c=0;c<=j;++c){
      Multiindex m1({c,j-c}), m2({0,c,j-c});
      for(int i=0;i<2;++i)
        P(i,m1) = jet(i+1,m2);
    }
  // Coefficient of the Birkhoff normal form should not vanish - see (*@\cite{SiegelMoser1971}@*)
  std::cout << "twist? " << 
    normalForms::computePlanarEllipticNormalForm(P)[1].real();
}
/* Output:
subset(N,r)? = true, N = [-2.9969206925386515e-15, 1.1389510487015546e-14]
twist? [15.406918966906115, 15.406919009018637]
*/
\end{lstlisting}

\section{Summary}
In this article we described a basic interface of the \capd library. The C++ source code of the \capd library consists of over 120\, 000 lines and thus it is clear that presenting all implemented features and details of algorithms in one article is impossible. We showed, however, that the library is a powerful tool for rigorous numerical analysis of dynamical systems by examining several non-trivial examples. 

We would like to mention, that the library provides support for integration of differential inclusions \cite{KapelaZgliczynski2009} and algorithms for rigorous integration of dissipative PDEs \cite{WilczakZgliczynski2020}. In the nearest future a $\mathcal C^1$ algorithm for PDEs, constrained $\mathcal C^0-\mathcal C^1$ algorithms for ODEs (that is taking into account constraints of the system, like Hamiltonians, measure preservation) and for delay differential equations should be added.

\bibliographystyle{elsarticle-num}
\addcontentsline{toc}{chapter}{Bibliography}
%uncomment next line to change bibliography name to references
\renewcommand{\bibname}{refs}
\section*{References}
\bibliography{refs}
\end{document}